\documentclass[12pt]{amsart}
\usepackage[alphabetic]{amsrefs}
\usepackage{mathdots, blkarray}
\usepackage{multicol}
\usepackage{graphicx}
\usepackage{amscd}
\usepackage{upgreek}
\usepackage{stmaryrd}
\usepackage{longtable}
\usepackage[T1]{fontenc}
\usepackage{latexsym, amsmath, amssymb, amsthm}
\usepackage[svgnames]{xcolor}
\usepackage{rsfso}
\usepackage[mathscr]{eucal}
\usepackage{mathtools}
\usepackage{mathptmx}
\usepackage{titletoc}
\usepackage{wrapfig}
\usepackage{float}
\usepackage{xypic}
\usepackage{microtype}
\usepackage{dsfont}
\usepackage{xcolor}
\usepackage{color}
\usepackage[colorlinks = true,
            linkcolor  = DarkBlue,
            urlcolor   = DarkRed,
            citecolor  = DarkGreen]{hyperref}
\usepackage{hyperref}
\allowdisplaybreaks	
\usepackage[centering, includeheadfoot, hmargin=1.0in, tmargin=0.5in, 
  bmargin=0.6in, headheight=6pt]{geometry}
\newtheorem{theorem}{Theorem}[section]
\newtheorem{prop}[theorem]{Proposition}
\newtheorem{defn}[theorem]{\rm\textsc{Definition}}
\newtheorem{lem}[theorem]{Lemma}
\newtheorem{coro}[theorem]{Corollary}
\newtheorem{conj}[theorem]{\rm \textsc{Conjecture}}
\newtheorem{thm}[theorem]{Theorem}

\newtheorem{rem}[theorem]{\rm\textsc{Remark}}
\newtheorem{exam}[theorem]{\rm\textsc{Example}}
\newcommand{\ideal}[1]{\ensuremath{\left\langle #1 \right\rangle}}

\newcommand{\oeq}[1]{\ensuremath{\overset{(\ref{#1})}{=}}}
\newcommand{\bslash}{\kern-0.1em\texttt{\scalebox{0.6}[1]{/}}\kern-0.15em \texttt{\scalebox{0.6}[1]{/}}}

\newcommand{\A}{\mathcal{A}} 
\newcommand{\NN}{\mathcal{N}} 
\newcommand{\C}{\mathbb{C}}

\newcommand{\N}{\mathbb{N}} 
 
\newcommand{\gl}{\mathfrak{gl}} 
\newcommand{\p}{\mathfrak{p}} 
\newcommand{\K}{\mathcal{K}} 
\newcommand{\KK}{\mathbb{K}} 
\newcommand{\LA}{\Longleftarrow} 
\newcommand{\RA}{\Longrightarrow} 
 
\newcommand{\ra}{\longrightarrow}

\newcommand{\hbo}{$\hfill\Diamond$} 
\newcommand{\dm}{\diamond}

\begin{document}
\title{Kupershmidt-Nijenhuis structures on pre-Malcev algebras} 
\def\shorttitle{Kupershmidt-Nijenhuis structures on pre-Malcev algebras}

\author{Yin Chen}
\address{School of Mathematics and Physics, Key Laboratory of ECOFM of 
Jiangxi Education Institute, Jinggangshan University,
Ji'an 343009, Jiangxi, China \& Department of Finance and Management Science, University of Saskatchewan, Saskatoon, SK, Canada, S7N 5A7}
\email{yin.chen@usask.ca}

\author{Liman Qin}
\address{School of Mathematics and Statistics, Northeast Normal University,
 Changchun, China}
\email{qinlm549@nenu.edu.cn}

\author{Shan Ren}
\address{School of Mathematics and Statistics, Northeast Normal University, Changchun, China}
\email{rens734@nenu.edu.cn}

\author{Runxuan Zhang}
\address{Department of Mathematics and Information Technology, Concordia University of Edmonton, Edmonton, AB, Canada, T5B 4E4}
\email{runxuan.zhang@concordia.ab.ca}

\begin{abstract}
We study Kupershmidt operators, Nijenhuis operators, and Kupershmidt-Nijenhuis structures on finite-dimensional pre-Malcev algebras over a field of characteristic zero. We construct several new families of complex pre-Malcev algebras that are not pre-Lie algebras in dimensions two, three and four. We use the compatibility of linear operators to establish  connections between Kupershmidt operators, Nijenhuis operators and Kupershmidt-Nijenhuis structures on pre-Malcev algebras. Moreover, we use a method from computational ideal theory to characterize the geometric structures of 
the varieties of Kupershmidt operators and Nijenhuis operators on a three-dimensional complex pre-Malcev algebra. 
\end{abstract}

\date{\today}
\maketitle \baselineskip=16.5pt

%%%%%%%%%%%%%%%%%%%%%%%%%%%Contents%%%%%%%%%%%%%%%%%%%%%%%%
%\textcolor{blue}{\tableofcontents{}}
\dottedcontents{section}[1.16cm]{}{1.8em}{5pt}
\dottedcontents{subsection}[2.00cm]{}{2.7em}{5pt}
%\dottedcontents{subsubsection}[2.86cm]{}{3.4em}{5pt}

%%%%%%%%%%%%%%%%%%%%%%%%%%%Sections%%%%%%%%%%%%%%%%%%%%%%%%
\section{Introduction}
\setcounter{equation}{0}
\renewcommand{\theequation}
{1.\arabic{equation}}
\setcounter{theorem}{0}
\renewcommand{\thetheorem}
{1.\arabic{theorem}}

\noindent Linear operators of finite-dimensional nonassociative algebras over a field have occupied a significant position in the structure theory of algebras. One of the main reasons for this is that properties of linear operators provide effective tools for describing the global structures of the underlying algebras. This perspective can be traced back to a classical theorem of Jacobson \cite{Jac55} on Lie derivations, which establishes that in the case of characteristic zero,  the existence of an invertible derivation implies the nilpotence of the Lie algebra. 
The theory of linear operators also has many substantial ramifications and applications in the areas of  functional analysis, 
mathematical physics, and differential geometry; see for example \cite{BGL25, CCZ21, GGZ21, HMNS23b} and \cite{Zha21}. The primary objective of this article is to explore
Kupershmidt operators, Nijenhuis operators, and Kupershmidt-Nijenhuis structures on finite-dimensional pre-Malcev algebras, a class of nonassociative algebras that includes all pre-Lie algebras as a special case. 

The concepts of Kupershmidt operators and Nijenhuis operators stem from mathematical physics and were initially introduced in the context of Lie algebras. Kupershmidt operators are closely connected to the classical Yang-Baxter equation (CYBE) and in fact, Kupershmidt \cite{Kup99} originally introduced these operators to understand $r$-matrices, which can be used to construct solutions of the CYBE in Lie algebras. Nijenhuis operators date back to the cohomology and deformation theory of Lie algebras in \cite{NR67,NR68} where they provided an efficient method for constructing new Lie algebras, and have influenced much recent research; for example, \cite{WSBL19,CZ26a} and \cite{CZ24b}. In 2018, \cite{HLS18} formulated the notion of Kupershmidt-Nijenhuis structures for Lie algebras that combines Kupershmidt operators and Nijenhuis operators, motivated by the desire to generalize $r$-matrix-Nijenhuis structures associated with the CYBE. Recently, Kupershmidt-Nijenhuis structures have been extended to other nonassociative algebras such as alternative algebras, pre-Lie algebras, and Malcev algebras; see \cite{Sun21, LW22} and \cite{Mab22}.

The notion of Malcev algebras (also called Moufang-Lie algebras) was introduced by Anatoly Malcev in the 1950s as a generalization of Lie algebras. They play an analogous role in the geometry of Moufang loops: the tangent space at the identity of a smooth Moufang loop forms a Malcev algebra, just as the tangent space at the identity of a Lie group forms a Lie algebra;
see for example, \cite{Mal55, Kuz71, Ker79}, and \cite{GRSS21}. 
All associative algebras, the octonions, and more generally all alternative algebras are Malcev-admissible; that is, the commutator $[x,y]:=xy-yx$ gives a Malcev algebra structure. Furthermore, all Lie-admissible algebras, including pre-Lie algebras and Okubo algebras, provide further examples of Malcev-admissible algebras.

As a generalization of pre-Lie algebras, pre-Malcev algebras are also Malcev-admissible and were introduced by \cite{Mad17}
to understand the splitting procedure occurred in \cite{BBGN13} to some anti-commutative and Malcev's identities.
In comparison with the cases of Lie algebras and other nonassociative algebras, the study of Kupershmidt-Nijenhuis structures on pre-Malcev algebras remains relatively new and not yet well understood. Only a few results are known; see for example, \cite{HMNS21} and \cite[Section 4]{RZ24}. This article aims to study Kupershmidt operators, Nijenhuis operators,
Nijenhuis pairs, and Kupershmidt-Nijenhuis structures on finite-dimensional pre-Malcev algebras over a field $\KK$ of characteristic zero. To better understand connections between these objects, we will adopt and extend an approach 
from computational ideal theory that has already been shown to be effective in our previous work
 \cite{CZ23,CRSZ26a} and \cite{CZ26b}.

We observe that many concepts of linear operators on nonassociative algebras, such as Rota-Baxter operators, Kupershmidt operators, and Nijenhuis operators, were initially defined for associative algebras or Lie algebras and subsequently generalized to broader classes of nonassociative algebras; see \cite{GP25,HMNS23a, CHM24, ZC25} and \cite{CRSZ26b}. This leads us to recall or reformulate some definitions related to Kupershmidt-Nijenhuis structures for general nonassociative algebras. 
This general framework provides a unified perspective for studying Kupershmidt-Nijenhuis structures on various types of nonassociative algebras. Moreover, we shall construct a number of new examples of low-dimensional pre-Malcev algebras
over the complex field $\C$ and present explicit descriptions to demonstrate how our approach offers a geometric understanding of the varieties of Kupershmidt operators and Nijenhuis operators on pre-Malcev algebras. These examples not only enrich the existing family of pre-Malcev algebras but also illustrate the effectiveness of our method in exposing deeper structural properties through a geometric viewpoint.

We organize this article as follows. Section \ref{sec2} consists of two subsections. In Subsection \ref{subsec2.1}, we recap and generalize the concepts of Kupershmidt operators, Nijenhuis operators, Nijenhuis pairs, and Kupershmidt-Nijenhuis structures for arbitrary finite-dimensional nonassociative algebras over $\KK$. Subsection \ref{subsec2.2} specializes in 
finite-dimensional pre-Malcev algebras. We explore fundamental properties of Nijenhuis operators (Propositions \ref{prop2.3} and \ref{prop2.5}), and in addition, construct a new family $A(a,b)$ of two-dimensional complex pre-Malcev algebras that are not pre-Lie algebras, demonstrating that the notion of pre-Malcev algebras is a 
proper generalization of pre-Lie algebras even in the two-dimensional case; see Example \ref{exam2.2}. We provide three examples to illustrate how to understand the structures of the varieties of Kupershmidt operators, Nijenhuis operators, and Nijenhuis pairs on $A(1,0)$ in terms of algebraic geometry; see Examples \ref{exam2.4}, \ref{exam2.6}, and \ref{exam2.7}.

In Section \ref{sec3}, we investigate the compatibility of Kupershmidt operators and demonstrate how this compatibility can be used to establish connections between Nijenhuis operators and Kupershmidt operators on pre-Malcev algebras; see 
Proposition \ref{prop3.2} and Theorem \ref{thm3.3}. Moreover, we provide a method using Kupershmidt operators to construct new pre-Malcev algebra structures on finite-dimensional bimodules in Proposition \ref{prop3.4}. 
This construction is of crucial importance to understanding Kupershmidt-Nijenhuis structures on pre-Malcev algebras; see the first statement of Theorem \ref{thm3.6}. More significantly, if a Kupershmidt operator $T$ and a Nijenhuis operator $N$, together with a linear map $S$, form a  Kupershmidt-Nijenhuis structure on a pre-Malcev algebra, then Theorem \ref{thm3.6} shows their fundamental properties and compatibility relationships.

 Section \ref{sec4} contains more examples and explicit computations. We propose a conjecture on 
 two-dimensional pre-Malcev algebras over $\C$, giving the general shape of a two-dimensional complex pre-Malcev algebras that is not a pre-Lie algebra. We also construct several new families of three-dimensional and four-dimensional pre-Malcev algebras. Furthermore, we completely characterize the geometric structure of the variety of
 Kupershmidt operators on a three-dimensional complex pre-Malcev algebra with respect to the adjoint representation, showing that it is a one-dimensional affine variety consisting of two one-dimensional irreducible components; see Corollary \ref{coro4.11}. Applying a similar method, we are able to compute the variety of Nijenhuis operators of the three-dimensional pre-Malcev algebra without going into the detailed proofs; see Theorem \ref{thm2}.

Throughout this article, we assume that $\KK$ denotes any field of characteristic zero and $\C$ denotes the complex field.
All algebras, vector spaces, and bimodules are finite-dimensional. We write $A$ for a finite-dimensional pre-Malcev algebra over $\KK$ (except for Subsection \ref{subsec2.1}) and $V$ for a finite-dimensional $A$-bimodule unless otherwise specified. We write $A_V$ for the action of $A$ on $V$ in which $\ell_x$ and $r_x$ denote the left and right actions of $x\in A$ on $V$ respectively.  We use $B\oplus C$ to denote the block diagonal matrix $\begin{pmatrix}
    B  & 0   \\
    0  &  C
\end{pmatrix}$
 if $B$ and $C$ are two matrices (or linear maps).

\section{Nijenhuis Operators and Nijenhuis Pairs}\label{sec2}
\setcounter{equation}{0}
\renewcommand{\theequation}
{2.\arabic{equation}}
\setcounter{theorem}{0}
\renewcommand{\thetheorem}
{2.\arabic{theorem}}

\noindent  In this section, we first recap some fundamental concepts and facts on Kupershmidt operators, Nijenhuis operators, Nijenhuis pairs, and Kupershmidt-Nijenhuis structures for a finite-dimensional nonassociative algebra over a field $\KK$. In the second part of this section, we focus on pre-Malcev algebras and present a new family $A(a,b)$ of two-dimensional pre-Malcev algebras over $\C$. In particular, we compute the varieties of Kupershmidt operators, Nijenhuis operators, and Nijenhuis pairs on $A(1,0)$. 

\subsection{Some linear operators on nonassociative algebras} \label{subsec2.1}
Throughout this preliminary subsection, we assume that $A$ denotes a finite-dimensional nonassociative algebra over  $\KK$ and $V$ is a finite-dimensional $A$-bimodule. 

We say that a linear map $T:V\ra A$ is a \textit{Kupershmidt operator} on $A_V$ if
\begin{equation}
\label{eq2.1}
T(u)T(v)=T\left(\ell_{T(u)}(v)+r_{T(v)}(u)\right)
\end{equation}
for all $u,v\in V$. A linear map $N:A\ra A$ is called  a \textit{Nijenhuis operator} on $A$ if 
\begin{equation}
\label{eq2.2}
N(x)N(y)=N(N(x)y+xN(y)-N(xy))
\end{equation}
for all $x,y\in A$. We write $\K(A_V)$ for the set of all Kupershmidt operators on $A_V$ and denote by $\NN(A)$ the 
set of all Nijenhuis operators on $A$. Usually, describing the structures of $\K(A_V)$ and $\NN(A)$ is a challenging 
task because they are in general not vector spaces. Thus a point of view from algebraic geometry, regarding them as affine varieties, might be a better way to understand their structures.  We say that $T_1,T_2\in \K(A_V)$  are \textit{compatible}
if $T_1+T_2\in \K(A_V)$.

Assume that $N$ is a Nijenhuis operator on $A$ and $S:V\ra V$ is a linear map. We say that $(N,S)$ 
is a \textit{Nijenhuis pair} on $A_V$ provided that
\begin{eqnarray}
\ell_{N(x)} S & = & S\ell_{N(x)}+S\ell_x S-S^2\ell_x \label{eq2.3}\\
r_{N(x)} S & = & Sr_{N(x)}+Sr_x S-S^2r_x \label{eq2.4}
\end{eqnarray}
for all $x\in A$.

Nijenhuis operators and Nijenhuis pairs naturally occur in the study of infinitesimal deformations of $A_V$. Consider a special subclass $\A$ of nonassociative algebras over $\KK$, for example, $\A$ could be the family of all Lie algebras or the family of all pre-Lie algebras. Suppose that $A\in\A$ with an $A$-bimodule $V$. Let $\upomega:A\times A\ra A$ denote a bilinear map and 
let $\upalpha$ and $\upbeta$ be two linear maps from $A$ to $\gl(V)$. Define the following 
$t$-parameterized family of  multiplications and linear maps:
\begin{eqnarray}
x\cdot_t y & := & xy+t\upomega(x,y) \\
\ell_x^{\upalpha,t} &:= & \ell_x+t\upalpha_x\\
r_x^{\upbeta,t} &:= & r_x+t\upbeta_x
\end{eqnarray}
where $x,y\in A$.  We say that $(\upomega,\upalpha,\upbeta)$ generates an \textit{infinitesimal deformation} of $A_V$ provided that $A_t:=(A, \cdot_t)\in\A$ and  together with $\ell_x^{\upalpha,t}$ and $r_x^{\upbeta,t}$, $V$ becomes an $A_t$-bimodule.  Infinitesimal deformations provide an effective method to construct new objects in $\A$; see for example, \cite[Section 3]{CZ24b}.

An important topic in the theory of infinitesimal deformations is to decide when two infinitesimal deformations
$(A, \upomega,\upalpha,\upbeta)$ and $(A, \upomega',\upalpha',\upbeta')$ are equivalent. 
We say that $(A, \upomega,\upalpha,\upbeta)$ and $(A, \upomega',\upalpha',\upbeta')$ are \textit{equivalent} if the two algebras $A_t$ are isomorphic and the two corresponding $A_t$-bimodule structures on $V$ are isomorphic as well. An infinitesimal deformation $(A, \upomega,\upalpha,\upbeta)$ is called \textit{trivial} if it is equivalent to $A_V$.
More precisely, $(A, \upomega,\upalpha,\upbeta)$ is trivial if there are two linear maps $N:A\ra A$ and $S:V\ra V$ such that
\begin{eqnarray}
\upomega(x,y)& = & N(x)y+xN(y)-N(xy) \\
N(\upomega(x,y))& = & N(x)N(y) \\
\upalpha_x&=&\ell_{N(x)}+\ell_x S-S\ell_x\\
S\upalpha_x&=&\ell_{N(x)}S\\
\upbeta_x&=&r_{N(x)}+r_x S-Sr_x\\
S\upbeta_x&=&r_{N(x)}S
\end{eqnarray}
for all $x,y\in A$. Apparently,  (\ref{eq2.2}) can be obtained from the first two relations above; the two relations in the middle produce (\ref{eq2.3}), and (\ref{eq2.4}) can be seen by combining the last two relations. 

Moreover, given  a Kupershmidt operator $T$ on $A_V$ and a Nijenhuis pair $(N,S)$, we define two multiplications $\dm^{T}$ and $\dm^T_S$ on $V$ by
\begin{eqnarray}
u\dm^{T}v&:=&\ell_{T(u)}(v)+r_{T(v)}(u) \label{eq2.14}\\
u\dm^{T}_Sv&:=& S(u)\dm^T v+u\dm^T S(v)-S(u\dm^T v)  \label{eq2.15}
\end{eqnarray}
for all $u,v\in V$. We shall see that $TS$ is also a Kupershmidt operator on $A_V$; see Theorem  \ref{thm2} below.

\begin{defn}{\rm
Let $T:V\ra A$ be a Kupershmidt operator of $A$ and $(N,S)$ be a Nijenhuis pair of $A_V$.
We say that $(T,N,S)$ is a {\it Kupershmidt-Nijenhuis structure} on $A$ provided that
$TS=NT$ and 
\begin{equation}
\label{eq2.16}
u\dm^{TS}v=u\dm^T_S v
\end{equation}
for all $u,v\in V$.
}\end{defn}

This concept was first introduced by \cite[Definition 3.3]{HLS18} for Lie algebras and recently, it was extended to Malcev algebras in \cite[Definition 4.4]{Mab22}.

\subsection{Nijenhuis operators and Nijenhuis pairs} \label{subsec2.2}

This subsection specializes in pre-Malcev algebras. Recall that a \textit{pre-Malcev algebra} over $\KK$ is a nonassociative algebra $A$ over $\KK$ whose multiplication satisfies 
$P_M(x,y,z,t)=0$, where
\begin{equation*}
\label{ }
P_M(x,y,z,t):=(yz-zy)(xt)+((xy-yx)z-z(xy-yx))t+y((xz-zx)t)-x(y(zt))+z(x(yt))
\end{equation*}
 for all $x,y,z,t\in A$; see \cite[Definition 4]{Mad17}. The commutator $[x,y]:=xy-yx$ defines a Malcev algebra structure $[A]$ on the underlying space of $A$. This means that pre-Malcev algebras are Malcev admissible. 

Note that Malcev algebras generalize the notion of Lie algebras, pre-Malcev algebras 
can be viewed as a generalization of pre-Lie algebras; see \cite[Remark 8]{Mad17}. 
We say that a pre-Malcev algebra is \textit{nontrivial} if it is not a pre-Lie algebra. 
There are only a few known nontrivial pre-Malcev algebras in dimensions 4 and 5; see \cite[Example 4.1]{RZ24} and \cite[Examples 3.1 and 3.2]{HMNS21}. Except for these, it seems no other nontrivial examples have been found. In the following example, we construct a family of two-dimensional nontrivial pre-Malcev algebras over $\C$.

\begin{exam}\label{exam2.2}
{\rm
Let $A(a,b)$ be a two-dimensional nonassociative algebra over $\C$ spanned by $\{e_1,e_2\}$ defined by the following  
multiplications:  
$$e_1e_1=ae_1+4abe_2,~ e_1e_2=-ae_2,~ e_2e_1=-be_1,~e_2e_2=e_1+be_2$$
where $a \in \C^\times$ and $b\in\C$. A direct verification shows that $A(a,b)$ is a pre-Malcev algebra but is not a pre-Lie algebra. In particular,
taking $a=1$ and $b=0$ gives the following pre-Malcev algebra $A(1,0)$:
$$e_1e_1=e_1,~ e_1e_2=-e_2,~ e_2e_1=0,~e_2e_2=e_1$$
which might be the simplest nontrivial pre-Malcev algebra over $\C$.
\hbo}\end{exam}

More examples of pre-Malcev algebras can be constructed from Malcev algebras and their Rota-Baxter operators (i.e.,  Kupershmidt operators associated to adjoint representations); see for example, \cite[Proposition 4.2]{RZ24}. 
Precisely, suppose that $(M,[-,-])$ is a finite-dimensional Malcev algebra over $k$ and $\upvarphi:M\ra M$ denotes a 
Rota-Baxter operator on $M$. Then $$x\cdot^{\upvarphi} y:=[\upvarphi(x),y]$$
defines a pre-Malcev algebra structure on $M$, denoted by $M_{\upvarphi}$. Furthermore, together with $\upvarphi$, 
Nijenhuis operators on $M$ commutating with $\upvarphi$ also are Nijenhuis operators on $M_{\upvarphi}$.

\begin{prop}\label{prop2.3}
Let $(M,[-,-])$ be a finite-dimensional Malcev algebra over $k$ and $\upvarphi$ be a Rota-Baxter operator on $M$.
Then $\{N\in\NN(M)\mid \upvarphi N=N\upvarphi\}$ is contained in $\NN(M_{\upvarphi})$.
\end{prop}

\begin{proof}
 For all $x, y \in M$ and $N\in\NN(M)$ with $\upvarphi N=N\upvarphi$, we have
\begin{eqnarray*}
N(x)\cdot^{\upvarphi}N(y)&=&[\upvarphi(N(x)), N(y)]=[N(\upvarphi(x)), N(y)]\\
&=&N\left([N(\upvarphi(x)), y]+ [\upvarphi(x), N(y)]- N([\upvarphi(x), y])\right)\\
&=&N\left([\upvarphi(N(x)), y]+ [\upvarphi(x), N(y)]- N([\upvarphi(x), y])\right)\\
&=&N\left(N(x) \cdot^{\upvarphi}y+x \cdot^{\upvarphi}N(y)-N(x \cdot^{\upvarphi}y)\right),
\end{eqnarray*}
which means that $N$ is a Nijenhuis operator on $M_{\upvarphi}$.
\end{proof}

Consider the pre-Malcev algebra $A(1,0)$ that appeared in Example \ref{exam2.2} and we would like to
compute the variety $\NN(A(1,0))$ of all Nijenhuis operators on $A(1,0)$ in terms of algebraic geometry. Here we do not give all the details because our approach will be explained meticulously in Section \ref{sec4}. 

\begin{exam}\label{exam2.4}
{\rm
Suppose that $N=\begin{pmatrix}
     x_{11} & x_{12}   \\
     x_{21} & x_{22} 
\end{pmatrix}\in \NN(A(1,0))$ denotes a generic element, and the action of $N$ on $A(1,0)$ is given by
$$N(e_1)=x_{11}e_1+ x_{12}e_2 \textrm{ and } N(e_2)=x_{21}e_1+ x_{22}e_2.$$ 
Together with the generating relations in $A(1,0)$, substituting this action into (\ref{eq2.2}) where $(x,y)$ runs over
$\{(e_1,e_1),(e_1,e_2),(e_2,e_1), (e_2,e_2)\}$, we obtain the following polynomial equations
\begin{eqnarray*}
x_{12}^2 + 2x_{12}x_{21} &=&0\\
x_{11}x_{12} - x_{12}x_{22}&=&0\\
x_{11}^2 - 2x_{1}x_{22} + x_{12}x_{21} + 2x_{21}^2 + x_{22}^2&=&0.
\end{eqnarray*}
The ideal $J$ generated by the corresponding three polynomials above is a radical ideal of the polynomial ring $\C[ x_{11}, x_{12},  x_{21},x_{22}]$. Hence,
$$\NN(A(1,0))=V(J).$$
Namely, $\NN(A(1,0))$ is the affine subvariety of $\C^4$ defined by $J$. Note that $J$ is not a prime ideal. Thus $\NN(A(1,0))$ is reducible.  In fact, using the classical Gr\"oebner basis method in computational ideal theory (see for example, \cite[Section 2]{CZ23}), we obtain that
$$\NN(A(1,0))=V(J)=V(\p_1)\cup V(\p_2)$$
where $\p_1=\ideal{x_{11} - x_{22}, x_{12} + 2x_{21}}$ and $\p_2=\ideal{x_{11}^2 - 2x_{11}x_{22} + 2x_{21}^2 + x_{22}^2, x_{12}}$ are two prime ideals. As direct consequences, a Nijenhuis operator $N_1$ in $V(\p_1)$ is of the following form
$$N_1=\begin{pmatrix}
     a& -2b   \\
     b & a
\end{pmatrix}$$ 
where $a,b\in\C$; moreover, a Nijenhuis operator $N_2$ in $V(\p_2)$ has the form
$$N_2=\begin{pmatrix}
     a& 0   \\
     b & c
\end{pmatrix}$$ 
where $a^2 - 2ac + 2b^2 + c^2=0$ and $a,b,c\in\C$. Therefore, we see that each Nijenhuis operator $N$ on $A(1,0)$ has  the form  either $N_1$ or $N_2$.
\hbo}\end{exam}

To expose connections between Nijenhuis operators and Nijenhuis pairs on pre-Malcev algebras, we need to recall
the definition of pre-Malcev bimodules. 

Suppose that $A$ is a finite-dimensional pre-Malcev algebra over $\KK$. A finite-dimensional vector space $V$ over $\KK$ is called an \textit{$A$-bimodule} if there are two linear maps $\ell$ and $r$ both from $A$ to $\gl(V)$ such that
\begin{eqnarray*}
r_x r_y r_z-r_x r_y\ell_z-r_x\ell_y r_z+r_x\ell_y\ell_z-r_{z(yx)}+\ell_y r_{zx}+\ell_{zy}r_x-\ell_{yz}r_x-\ell_z r_x\ell_y+\ell_z r_x r_y&=& 0\\
r_{x}r_{y}\ell_{z}-r_{x}r_{y}r_{z}-r_{x}\ell_{y}\ell_{z}+r_{x}\ell_{y}r_{z}-\ell_{z}r_{yx}+\ell_{y}\ell_{z}r_{x}+r_{zx}r_{y}
-r_{zx}\ell_{y}-r_{(yz)x}+r_{(zy)x}&=& 0\\
r_{x}\ell_{yz}-r_{x}\ell_{zy}-r_{x}r_{yz}+r_{x}r_{zy}-\ell_{y}\ell_{z}r_{x}+r_{y(zx)}+r_{yx}\ell_{z}
-r_{yx}r_{z}-\ell_{z}r_{x}r_{y}+\ell_{z}r_{x}\ell_{y}&=& 0\\
\ell_{(xy)z}-\ell_{(yx)z}-\ell_{z(xy)}+\ell_{z(yx)}-\ell_{x}\ell_{y}\ell_{z}+\ell_{z}\ell_{x}\ell_{y}+\ell_{yz}\ell_{x}-\ell_{zy}\ell_{x}-\ell_{y}\ell_{zx}+\ell_{y}\ell_{xz}&=& 0
\end{eqnarray*}
where $x,y,z\in A$ and $\ell_x$ (or $r_x$) denotes the image of $x$ under $\ell$ (or $r$). Verifying these conditions is usually complicated and tedious. However, \cite[Section 4.2]{RZ24} asserts that 
it is equivalent to saying that the multiplication defined by 
\begin{equation}\label{eq2.17}
(x,u)(y,v):=(xy,\ell_x(v)+r_y(u))
\end{equation}
for $x,y\in A$ and $u,v\in V$, produces a pre-Malcev algebra structure on $A\oplus V$, which is called the \textit{semi-direct product} of $A$ and $V$
and denoted as $A\ltimes V$.

\begin{prop} \label{prop2.5}
Let $(N,S)$ be a Nijenhuis pair on a pre-Malcev algebra $A$ with respect to an $A$-bimodule $V$. Then $N\oplus S\in\NN(A\ltimes V)$.
\end{prop}

\begin{proof}
Suppose that $x,y\in A$ and $u,v\in V$ denote arbitrary elements. It follows from (\ref{eq2.17}) that 
$$(N\oplus S)(x,u)(N\oplus S)(y,v)=(N(x), S(u))(N(y),S(v))=(N(x)N(y), \ell_{N(x)}(S(v))+r_{N(y)}(S(u))).$$
Moreover, 
\begin{eqnarray*}
&&\hspace{-0.4cm}(N\oplus S)((N\oplus S)(x, u)(y, v)+(x,u)(N\oplus S)(y,v)-(N\oplus S)((x,u)(y,v)))\\
&\hspace{-0.6cm}=&\hspace{-0.4cm}(N\oplus S)((N(x), S(u))(y, v)+(x,u)(N(y),S(v))-(N\oplus S)(xy,\ell_{x}(v)+r_{y}(u))\\
&\hspace{-0.6cm}=&\hspace{-0.4cm}(N\oplus S)((N(x)y,\ell_{N(x)}(v)+r_y(S(u)))+(xN(y),\ell_{x}(S(v))+r_{N(y)}(u))-(N(xy),S(\ell_{x}(v)+r_{y}(u))))\\
&\hspace{-0.6cm}=&\hspace{-0.4cm} (\Delta_1, \Delta_2)
\end{eqnarray*}
where $\Delta_1:=N(N(x)y+xN(y)-N(xy))$ and $\Delta_2:=S(\ell_{N(x)}(v)+\ell_{x}(S(v))-S(\ell_{x}(v)))+ S(r_{N(y)}(u)+r_{y}(S(u))-S(r_{y}(u)))$. Since $(N,S)$ is a Nijenhuis pair, we see that 
$$(N\oplus S)(x,u)(N\oplus S)(y,v)=(N(x)N(y), \ell_{N(x)}(S(v))+r_{N(y)}(S(u)))=(\Delta_1, \Delta_2).$$
Hence, $N\oplus S$ is a Nijenhuis operator on $A\ltimes V$.
\end{proof}

We are interested in Nijenhuis pairs of pre-Malcev algebras with respect to adjoint bimodules because 
in this case, Kupershmidt operators coincide with the notion of Rota-Baxter operators and the latter have been extensively
studied for various nonassociative algebras; see for example, \cite{AB08, Guo12, GK08, GLS21} and \cite{TZS14}.

Suppose that $V=A$ denotes the adjoint bimodule of a pre-Malcev algebra $A$.  A Kupershmidt operator on $A$ is a linear map $T:A\ra A$ satisfying 
\begin{equation}
\label{eq2.18}
T(x)T(y)=T\left(T(x)y+xT(y)\right)
\end{equation}
for all $x,y\in A$. Moreover, a Nijenhuis pair on $A$ is a pair  $(N,S)$ satisfying 
\begin{eqnarray}
N(x)S(y) & = & S(N(x)y+xS(y)-S(xy)) \label{eq2.19}\\
S(y)N(x) & = & S(yN(x)+S(y)x-S(yx))  \label{eq2.20}
\end{eqnarray}
where $x,y\in A$, $N\in\NN(A)$ and $S:A\ra A$ is a linear map.  

In the next section, we shall see that Kupershmidt operators and Nijenhuis pairs play a fundamental role in understanding Kupershmidt-Nijenhuis structures on pre-Malcev algebras. Thus we close this section by providing two examples that compute the varieties of Kupershmidt operators and Nijenhuis pairs of $A(1,0)$ in Example \ref{exam2.2} with respect to the adjoint bimodule.

\begin{exam}\label{exam2.6}
{\rm
Let $T=\begin{pmatrix}
     x_{11} & x_{12}   \\
     x_{21} & x_{22} 
\end{pmatrix}\in \K(A(1,0))$ be an arbitrary Kupershmidt operator on $A(1,0)$ associated to the adjoint representation. 
Suppose that the action of $T$ on $A(1,0)$ is given by
$$T(e_1)=x_{11}e_1+ x_{12}e_2 \textrm{ and } T(e_2)=x_{21}e_1+ x_{22}e_2.$$ 
Substituting this action into (\ref{eq2.18}) and making $(x,y)$ run over
$\{(e_1,e_1),(e_1,e_2),(e_2,e_1), (e_2,e_2)\}$, we obtain six polynomials equations
\begin{eqnarray*}
2x_{21}^2 + x_{22}^2=0 && x_{12}x_{21} + x_{22}^2=0\\
 x_{12}^2 + 2x_{22}^2=0&& 3x_{11}x_{21} + 2x_{12}x_{22} - x_{21}x_{22}=0\\
 x_{11}x_{12} + x_{12}x_{22}=0&&x_{11}^2 - x_{22}^2=0.
\end{eqnarray*}
The ideal $J$ generated by the corresponding six polynomials is not a radical ideal of  $\C[ x_{11}, x_{12},  x_{21},x_{22}]$. Hence, $\K(A(1,0))=V(\sqrt{J})$, where
$$\sqrt{J}=\ideal{x_{11} + x_{22}, x_{12} - 2x_{21}, 2x_{21}^2 + x_{22}^2}$$
is actually a prime ideal. Therefore, $\K(A(1,0))$ is an irreducible affine variety of Krull dimension 1 whose elements have the following form 
$$T=\begin{pmatrix}
    a& 2b  \\
     b& -a
\end{pmatrix}$$
where $a,b\in\C$ and $a^2+2b^2=0.$
\hbo}\end{exam}

A similar method can be applied to compute all Nijenhuis pairs $(N_1,S)$ on $A(1,0)$, where
$$N_1=\begin{pmatrix}
     a& -2b   \\
     b & a
\end{pmatrix}$$ 
denotes a Nijenhuis operator on $A(1,0)$ occurred in Example \ref{exam2.4}.

\begin{exam}\label{exam2.7}
{\rm
Suppose that $S=\begin{pmatrix}
     x_{11} & x_{12}   \\
     x_{21} & x_{22} 
\end{pmatrix}$ denotes a linear map on $A(1,0)$ with the following action on $A(1,0)$ defined by: $$S(e_1)=x_{11}e_1+ x_{12}e_2 \textrm{ and } S(e_2)=x_{21}e_1+ x_{22}e_2.$$ To determine all $S$ such that $N_1$ and $S$ form a  Nijenhuis pair, it suffices to substitute the actions of
$S$ and $N_1$ on $A(1,0)$ to (\ref{eq2.18}) and (\ref{eq2.19}) and we obtain a finite set of polynomial equations in $x_{11}, x_{12}, x_{21}, x_{22}$. The ideal $J$ generated by the corresponding polynomials 
determines an affine variety that has two irreducible components $V(\p_1)$ and $V(\p_2)$, where
$$\p_1=\ideal{x_{11}-x_{22},x_{12},x_{21}}\textrm{ and }\p_2=\ideal{x_{11}-a,x_{12}+2b,x_{21}-b,x_{22}-a}.$$
Hence, $S$ either equals $N_1$ or has the form $\begin{pmatrix}
    c & 0  \\
     0 & c
\end{pmatrix}$, where $c\in\C$.
\hbo}\end{exam}

\section{Kupershmidt-Nijenhuis Structures} \label{sec3}
\setcounter{equation}{0}
\renewcommand{\theequation}
{3.\arabic{equation}}
\setcounter{theorem}{0}
\renewcommand{\thetheorem}
{3.\arabic{theorem}}

\noindent  This section examines Kupershmidt–Nijenhuis structures on a finite-dimensional pre-Malcev algebra $A$ (with respect to a bimodule $V$) over $\KK$ and explores some connections between Kupershmidt operators and Nijenhuis operators. A key element in establishing these  connections is the compatibility between two operators. 

We begin with the following result.

\begin{lem}\label{lem3.1}
 Two Kupershmidt operators $T_1$ and $T_2$ on $A_V$ are compatible  if and only if 
 \begin{equation}\label{eq3.1}
T_{1}(u)T_{2}(v)+T_{2}(u)T_{1}(v)
=T_{1}\left(\ell_{T_{2}(u)}(v)+r_{T_{2}(v)}(u)\right)+T_{2}\left(\ell_{T_{1}(u)}(v)+r_{T_{1}(v)}(u)\right)
\end{equation}
for all $u,v\in V.$
\end{lem}

\begin{proof}
$(\RA)$ Assume that $T_{1}$ and $T_{2}$ are compatible, i.e., $T_{1}+T_{2}\in\K(A_V)$. Thus
\begin{eqnarray*}
(T_{1}+T_{2})(u)(T_{1}+T_{2})(v)&=&(T_{1}+T_{2})(\ell_{(T_{1}+T_{2})(u)}(v)+r_{(T_{1}+T_{2})(v)}(u))\\
&=&(T_{1}+T_{2})(\ell_{T_{1}(u)}(v)+\ell_{T_{2}(u)}(v)+r_{T_{1}(v)}(u)+r_{T_{2}(v)}(u))\\
&=&T_{1}(\ell_{T_{1}(u)}(v)+r_{T_{1}(v)}(u))+T_{2}(\ell_{T_{2}(u)}(v)+r_{T_{2}(v)}(u))+\\
&&T_{1}(\ell_{T_{2}(u)}(v)+r_{T_{2}(v)}(u))+T_{2}(\ell_{T_{1}(u)}(v)+r_{T_{1}(v)}(u)).
\end{eqnarray*}
On the other hand, note that $(T_{1}+T_{2})(u)(T_{1}+T_{2})(v)=T_{1}(u)T_{1}(v)+T_{2}(u)T_{2}(v)+T_{1}(u)T_{2}(v)+T_{2}(u)T_{1}(v)$ and $T_1,T_2\in\K(A_V)$, it follows (\ref{eq2.1}) that
$$T_{1}(u)T_{2}(v)+T_{2}(u)T_{1}(v)
=T_{1}\left(\ell_{T_{2}(u)}(v)+r_{T_{2}(v)}(u)\right)+T_{2}\left(\ell_{T_{1}(u)}(v)+r_{T_{1}(v)}(u)\right).$$
$(\LA)$ Conversely, substituting $T=T_1$ and $T_2$ in (\ref{eq2.1}) and together with (\ref{eq3.1}), we see that
$$(T_{1}+T_{2})(u)(T_{1}+T_{2})(v)= (T_{1}+T_{2})(\ell_{(T_{1}+T_{2})(u)}(v)+r_{(T_{1}+T_{2})(v)}(u)).$$
This means that $T_1+T_2\in\K(A_V)$ and thus $T_{1}$ and $T_{2}$ are compatible.
\end{proof}

Invertible Nijenhuis operators can be used to construct new compatible Kupershmidt operators.

\begin{prop}\label{prop3.2}
Suppose that $T\in\K(A_V)$ and $N\in\NN(A)$. Then $NT\in\K(A_V)$ if and only if
\begin{equation}
\label{eq3.2}
N\left(NT(u)T(v)+T(u)NT(v)\right)=N\left(T(\ell_{NT(u)}(v)+r_{NT(v)}(u))+NT(\ell_{T(u)}(v)+r_{T(v)}(u))\right)
\end{equation}
for all $u,v\in V.$ In this case, if $N$ is invertible, then $T$ and $ NT $ are compatible.
\end{prop}

\begin{proof}
Let us prove the first statement. Since $N$ is a Nijenhuis operator on $A$, it follows from (\ref{eq2.2}) with $(x,y)=(T(u),T(v))$ that 
$$
NT(u)\cdot NT(v) = N\left(NT(u)T(v)+T(u)NT(v)-N(T(u)T(v))\right).
$$
The assumption that
$T\in \K(A_V)$ implies that $T(u)T(v)=T\left(\ell_{T(u)}(v)+r_{T(v)}(u)\right).$ Hence,
\begin{equation}
\label{eq3.3}
NT(u)\cdot NT(v) = N\left(NT(u)T(v)+T(u)NT(v)-NT(\ell_{T(u)}(v)+r_{T(v)}(u))\right).
\end{equation}
Combining (\ref{eq3.2}) and (\ref{eq3.3}), we observe that the condition (\ref{eq3.2}) holds if and only if
$$NT(u)\cdot NT(v) =NT(\ell_{NT(u)}(v)+r_{NT(v)}(u))$$
which is equivalent to saying that $NT$ is a Kupershmidt operator on $A_V$.

To see the second statement, multiplying $N^{-1}$ on both sides of (\ref{eq3.2}) implies that
$$NT(u)T(v)+T(u)NT(v)=T(\ell_{NT(u)}(v)+r_{NT(v)}(u))+NT(\ell_{T(u)}(v)+r_{T(v)}(u)).$$
By Lemma \ref{lem3.1}, it follows that $T$ and $ NT $ are compatible.
\end{proof}

The following result  establishes a strong connection between 
invertible and compatible Kupershmidt operators and Nijenhuis operators.

\begin{thm}\label{thm3.3}
Suppose that $T_1,T_2\in\K(A_V)$ both are invertible. Then $T_1$ and $T_2$ are compatible if 
and only if $T_1T_2^{-1}\in\NN(A)$.
\end{thm}

\begin{proof}
To begin with, let us define $N:=T_1T_2^{-1}$ and assume that $x,y\in A$ are two arbitrary elements.

$(\RA)$  We first note that there are two elements $u,v\in V$ such that $T_2(u)=x$ and $T_2(v)=y$. Since $NT_2=T_1\in\K(A_V)$, it follows that
$$N(T_{2}(u))N(T_{2}(v))=N\left(T_{2}\left(\ell_{N(T_{2}(u))}(v)+r_{N(T_{2}(v))}(u)\right)\right).$$
Note that $T_2\in\K(A_V)$ and it is compatible with $T_1$, thus it follows from Lemma \ref{lem3.1} that
\begin{eqnarray*}
&&N(T_{2}(u))\cdot T_{2}(v)+T_{2}(u)\cdot N(T_{2}(v))\\
&=&N\left(T_{2}\left(\ell_{T_{2}(u)}(v)+r_{T_{2}(v)}(u)\right)\right)+T_{2}\left(\ell_{N(T_{2}(u))}(v)+r_{N(T_{2}(v))}(u)\right)\\
&\oeq{eq2.1}&T_{2}\left(\ell_{N(T_{2}(u))}(v)+r_{N(T_{2}(v))}(u)\right)+N\left(T_{2}(u)\cdot T_{2}(v)\right).
\end{eqnarray*}
Letting $N$ act on both sides gives
$$N\left(N(T_{2}(u))\cdot T_{2}(v)+T_{2}(u)\cdot N(T_{2}(v))\right)=N\left(T_{2}\left(\ell_{N(T_{2}(u))}(v)+r_{N(T_{2}(v))}(u)\right)\right)+N^2\left(T_{2}(u)\cdot T_{2}(v)\right)$$
which equals $N(T_{2}(u))N(T_{2}(v))+N^{2}\left(T_{2}(u)T_{2}(v)\right)$. Hence,
\begin{equation*}
\label{ }
N(T_{2}(u))N(T_{2}(v))=N\left(N(T_{2}(u))\cdot T_{2}(v)+T_{2}(u)\cdot N(T_{2}(v))\right)-N^{2}\left(T_{2}(u)T_{2}(v)\right).
\end{equation*}
Namely, $N(x)N(y)=N(N(x)y+xN(y)-N(xy))$ and so $N\in\NN(A)$.

$(\LA)$  The converse statement can be seen by  Proposition \ref{prop3.2}. In fact,
note that $T_1$ and $T_2$ both are invertible, thus $N=T_1T_2^{-1}$ is also invertible. Since $NT_2=T_1$ is a 
Kupershmidt operator, it follows from the second statement in Proposition \ref{prop3.2} that
$T_2$ and $T_1$ are compatible. 
\end{proof}

We may use Kupershmidt operators to construct new pre-Malcev algebra structures.

\begin{prop}\label{prop3.4}
Let $V$ be a finite-dimensional $A$-bimodule and $T\in \K(A_V)$. Then $V$ endowed with the multiplication $u\dm^{T}v$ in $(\ref{eq2.14})$ is a pre-Malcev algebra.
\end{prop}

\begin{proof}
To see that $(V, u\dm^{T}v)$ is a pre-Malcev algebra, we need to verify $P_M(u,v,w,t)=0$ for all $u, v, w, t\in V$, where
\begin{eqnarray*}
P_M(u,v,w,t)&\hspace{-2mm}=(v\dm^{T} w)\dm^{T}(u\dm^{T} t)-(w\dm^{T} v)\dm^{T}(u\dm^{T} t)+((u\dm^{T} v)\dm^{T} w)\dm^{T} t-((v\dm^{T} u)\dm^{T} w)\dm^{T} t\\
&\hspace{-3cm} + (w\dm^{T}(v\dm^{T} u))\dm^{T} t-(w\dm^{T}(u\dm^{T} v))\dm^{T} t+v\dm^{T}((u\dm^{T} w)\dm^{T} t)\\
&\hspace{-3cm} -v\dm^{T}((w\dm^{T} u)\dm^{T} t+w\dm^{T}(u\dm^{T}(v\dm^{T} t))-u\dm^{T}(v\dm^{T}(w\dm^{T} t)).
\end{eqnarray*}
According to (\ref{eq2.14}) and expanding all terms on the right-hand side of $P_M(u,v,w,t)$, we see that
\begin{eqnarray*}
&&(v\dm^{T} w)\dm^{T}(u\dm^{T} t)-(w\dm^{T} v)\dm^{T}(u\dm^{T} t)\\
&=&\ell_{T(v)T(w)-T(w)T(v)}(\ell_{T(u)}(t)+r_{T(t)}(u))+r_{T(u)T(t)}(\ell_{T(v)}(w)+r_{T(w)}(v)-\ell_{T(w)}(v)-r_{T(v)}(w)).
\end{eqnarray*}
Similarly, 
\begin{eqnarray*}
&&((u\dm^{T} v)\dm^{T} w)\dm^{T} t-((v\dm^{T} u)\dm^{T}w)\dm^{T} t\\
&=&\ell_{T\left(\ell_{T(u)T(v)-T(v)T(u)}(w)+r_{T(w)}(\ell_{T(u)}(v)+r_{T(v)}(u)-\ell_{T(v)}(u)-r_{T(u)}(v)\right)}(t)+\\
&&r_{T(t)}\left(\ell_{T(u)T(v)-T(v)T(u)}(w)+r_{T(w)}(\ell_{T(u)}(v)+r_{T(v)}(u)-\ell_{T(v)}(u)-r_{T(u)}(v))\right),
\end{eqnarray*}
\begin{eqnarray*}
&&(w\dm^{T}(v\dm^{T} u))\dm^{T} t-(w\dm^{T}(u\dm^{T} v))\dm^{T} t\\
&=&\ell_{T\left(\ell_{T(w)}(\ell_{T(v)}(u)-\ell_{T(u)}(v)+r_{T(u)}(v)-r_{T(v)}(u))+r_{T(v)T(u)-T(u)T(v)}(w)\right)}(t)+\\
&&r_{T(t)}\left(\ell_{T(w)}(\ell_{T(v)}(u)-\ell_{T(u)}(v)+r_{T(u)}(v)-r_{T(v)}(u))+r_{T(v) T(u)-T(u)T(v)}(w)\right),
\end{eqnarray*}
\begin{eqnarray*}
&&v\dm^{T}((u\dm^{T} w)\dm^{T} t)-v\dm^{T}((w\dm^{T}u)\dm^{T} t)\\
&=&r_{T\left(\ell_{T(u)T(w)-T(w)T(u)}(t)+r_{T(t)}(\ell_{T(u)}(w)+r_{T(w)}(u)-\ell_{T(w)}(u)-r_{T(u)}(w))\right)}(v)+\\
&&\ell_{T(v)}\left(\ell_{T(u)T(w)-T(w)T(u)}(t)+r_{T(t)}(\ell_{T(u)}(w)+r_{T(w)}(u)-\ell_{T(w)}(u)-r_{T(u)}(w))\right),
\end{eqnarray*} and
\begin{eqnarray*}
&&\hspace{-2mm}w\dm^{T}(u\dm^{T}(v\dm^{T} t))-u\dm^{T}(v\dm^{T}(w\dm^{T}t))\\
&\hspace{-3mm}=&\hspace{-2mm}r_{T\left(\ell_{T(u)}(\ell_{T(v)}(t)+r_{T(t)}(v))+r_{T(v)T(t)}(u)\right)}(w)-r_{T\left(\ell_{T(v)}(\ell_{T(w)}(t)+r_{T(t)}(w))+r_{T(w)T(t)}(v)\right)}(u)+\\
&&\hspace{-2mm}\ell_{T(w)}(\ell_{T(u)}(\ell_{T(v)}(t)+r_{T(t)}(v))+r_{T(v)T(t)}(u))-\ell_{T(u)}(\ell_{T(v)}(\ell_{T(w)}(t)+r_{T(t)}(w))
+r_{T(w)T(t)}(v)).
\end{eqnarray*}
Substituting these equations together obtains
\begin{eqnarray*}
P_M(u,v,w,t)&=& (\ell_{T(v)T(w)}\ell_{T(u)}-\ell_{T(w)T(v)}\ell_{T(u)}+\ell_{(T(u)T(v))T(w)}-\ell_{(T(v) T(u)) T(w)}\\
&&+ \ell_{T(w)(T(v)T(u))}-\ell_{T(w)(T(u) T(v))}+ \ell_{T(v)}\ell_{T(u)T(w)}-\ell_{T(v)}\ell_{T(w)T(u)}\\
&&+ \ell_{T(w)}\ell_{T(u)}\ell_{T(v)}-\ell_{T(u)}\ell_{T(v)}\ell_{T(w)})(t)\\
&&\hspace{-3mm}+ (\ell_{T(v)T(w)}r_{T(t)}-\ell_{T(w)T(v)}r_{T(t)}-r_{T(t)}r_{T(w)}\ell_{T(v)}+r_{T(t)}r_{T(w)}r_{T(v)}\\
&&+ r_{T(t)}\ell_{T(w)}\ell_{T(v)}-r_{T(t)}\ell_{T(w)}r_{T(v)}-\ell_{T(v)}r_{T(t)}\ell_{T(w)}+\ell_{T(v)}r_{T(t)}r_{T(w)}\\
&&+ \ell_{T(w)}r_{T(v)T(t)}-r_{T(v)(T(w)T(t))})(u)\\
&&\hspace{-3mm} + (r_{T(u)T(t)}r_{ T(w)}-r_{T(u)T(t)}\ell_{ T(w)}+r_{T(t)}r_{T(w)}\ell_{ T(u)}-r_{T(t)}r_{T(w)}r_{T(u)}\\
&&- r_{T(t)}\ell_{T(w)}\ell_{T(u)}+r_{T(t)}\ell_{T(w)}r_{ T(u)}+r_{(T(u) T(w)) T(t)}-r_{(T(w) T(u))T(t)}\\
&& +\ell_{T(w)}\ell_{T(u)}r_{T(t)}-\ell_{T(u)}r_{T(w)T(t)})(v)\\
&&\hspace{-3mm}+ (r_{T(u)T(t)}\ell_{T(v)}-r_{T(u)T(t)}r_{T(v)}+r_{T(t)}\ell_{T(u)T(v)}-r_{T(t)}\ell_{T(v)T(u)}\\
&&+ r_{T(t)}r_{T(v)T(u)}-r_{T(t)}r_{T(u)T(v)}+\ell_{T(v)}r_{T(t)}\ell_{T(u)}-\ell_{T(v)}r_{T(t)}r_{T(u)}\\
&&+ r_{T(u) (T(v)T(t))} -\ell_{T(u)}\ell_{T(v)}r_{T(t)})(w).
\end{eqnarray*}
Since $V$ is a pre-Malcev bimodule, it follows that $P_M(u,v,w,t)=0$ for all $u, v, w, t\in V$.
Hence, $(V,\dm^{T})$ is a pre-Malcev algebra.
\end{proof}

\begin{rem}{\rm
Combining (\ref{eq2.1}),  (\ref{eq2.14}), and Proposition \ref{prop3.4}, we see that 
\begin{equation}
\label{eq3.4}
T(u\dm^T v)=T(u)T(v)
\end{equation}
for all $T\in\K(A_V)$ and $u,v\in V$. In other words, $T$ is a homomorphism of pre-Malcev algebras. 
\hbo}\end{rem}

The following theorem is the main result in this section that reveals some fundamental relationships between 
Kupershmidt operators and Nijenhuis operators in a Kupershmidt-Nijenhuis structure on pre-Malcev algebras.

\begin{thm}\label{thm3.6}
Let $(T,N,S)$ be a Kupershmidt-Nijenhuis structure on $A_V$. Suppose that $u$ and $v$ are two arbitrary elements in $V$. Then 
\begin{enumerate}
  \item $S\in\NN(V,\dm^{T})$;
  \item $NT\in\K(A_V)$;
  \item $T$ and $TS$ are compatible Kupershmidt operators on $A_V$;
  \item For $k,i\in\N$, define $T_k:=TS^k=N^kT$ and 
  \begin{equation*}
\label{ }
T(u)\cdot_{N^{i}}T(v):=N^{i}(T(u))\cdot T(v)+T(u)\cdot N^{i}(T(v))-N^{i}(T(u)\cdot T(v)).
\end{equation*}
Then $T_k\in\K(A_V)$ and $T_k(u\dm^T_{S^{k+i}} v)=T_{k}(u)\cdot_{N^i}T_{k}(v)$.
  \item Any two $T_k, T_j$ are compatible for $k,j\in\N$.
\end{enumerate}
\end{thm}

\begin{proof} 
(1) To show that $S\in\NN(V,\dm^{T})$, it suffices to verify that 
$$\Delta:=S(u)\dm^{T}S(v)-S(S(u)\dm^{T}v+u\dm^{T}S(v)-S(u\dm^{T}v))=0.$$
First of all, it follows from (\ref{eq2.14}) that $S(u)\dm^{T}S(v)=\ell_{TS(u)}(S(v))+r_{TS(v)}(S(u))$. Moreover, 
\begin{equation}
\label{eq3.5}
S(u)\dm^{T}v+u\dm^{T}S(v)-S(u\dm^{T}v)\oeq{eq2.15}u\diamond_{S}^{T}v\oeq{eq2.16}u\dm^{TS}v.
\end{equation}
Note that $u\dm^{TS}v=\ell_{TS(u)}(v)+r_{TS(v)}(u)$, thus
\begin{eqnarray*}
\Delta&=&\ell_{TS(u)}S(v)+r_{TS(v)}S(u)-S(u\dm^{TS}v) \\
&=&\ell_{TS(u)}S(v)+r_{TS(v)}S(u)- S(\ell_{TS(u)}(v)+r_{TS(v)}(u))\\
&=&\ell_{TS(u)}S(v)-S\ell_{TS(u)}(v)+r_{TS(v)}S(u)-Sr_{TS(v)}(u).
\end{eqnarray*}
To show that $\Delta=0$, we need more identities. Substituting $x=T(u)$ in (\ref{eq2.3}) obtains that 
$$\ell_{NT(u)} S(v) - S\ell_{NT(u)}(v)=S\ell_{T(u)} S(v)-S^2\ell_{T(u)}(v).$$
Similarly, we substitute $x=T(v)$ in (\ref{eq2.4}) and see that
$$r_{NT(v)} S(u) - Sr_{NT(v)}(u)=Sr_{T(v)} S(u)-S^2r_{T(v)}(u).$$
Since $TS=NT$, it follows that
\begin{eqnarray*}
\Delta&=&Sr_{T(v)} S(u)-S^2r_{T(v)}(u)+S\ell_{T(u)} S(v)-S^2\ell_{T(u)}(v) \\
&=&S\left(r_{T(v)} S(u)-Sr_{T(v)}(u)+\ell_{T(u)} S(v)-S\ell_{T(u)}(v)\right). 
\end{eqnarray*}
Thus, to show $\Delta=0$, it suffices to show that $r_{T(v)} S(u)-Sr_{T(v)}(u)+\ell_{T(u)} S(v)-S\ell_{T(u)}(v)=0$. In fact, this can be seen by (\ref{eq3.5}).  More precisely, note that $S(u)\dm^{T}v+u\dm^{T}S(v)-S(u\dm^{T}v)=\ell_{TS(u)}(v)+r_{T(v)}S(u)+\ell_{T(u)}S(v)+\ell_{TS(v)}(u)-S(\ell_{T(u)}(v)+r_{T(v)}(u))$. Thus (\ref{eq3.5}) can be read as
 \begin{equation}
\label{eq3.6}
S\ell_{T(u)}(v)+Sr_{T(v)}(u)=\ell_{T(u)}S(v)+r_{T(v)}S(u)
\end{equation}
which is equivalent to that $r_{T(v)} S(u)-Sr_{T(v)}(u)+\ell_{T(u)} S(v)-S\ell_{T(u)}(v)=0$. Therefore, $\Delta=S(0)=0$, as desired. 

(2) We need to prove that $NT(u)\cdot NT(v)=NT(\ell_{NT(u)}(v)+r_{NT(v)}(u)).$
Note that $NT(\ell_{NT(u)}(v)+r_{NT(v)}(u))=NT(u\dm^{NT}v)=NT(u\dm^{TS}v)=NT(u\dm_{S}^{T}v).$
Thus, it suffices to show that $$NT(u\dm_{S}^{T}v)=NT(u)\cdot NT(v).$$
In fact, 
\begin{eqnarray*}
NT(u\dm_{S}^{T}v)&\hspace{-3mm}=&\hspace{-3mm} NT\left(S(u)\dm^{T}v+u\dm^{T}S(v)-S(u\dm^{T}v)\right)\\
&\hspace{-3mm}=&\hspace{-3mm}NT\left(\ell_{TS(u)}(v)+r_{T(v)}S(u)+\ell_{T(u)}S(v)+r_{TS(v)}(u)-S(\ell_{T(u)}(v)+r_{T(v)}(u))\right)\\
&\hspace{-3mm}=&\hspace{-3mm}N\left(T(\ell_{TS(u)}(v)+r_{T(v)}S(u))+T(\ell_{T(u)}S(v)+r_{TS(v)}(u))-NT(\ell_{T(u)}(v)+r_{T(v)}(u))\right)\\
&\hspace{-3mm}\oeq{eq2.1}&\hspace{-3mm}N\left(T(S(u))T(v)+T(u)T( S(v))-N(T(u)T(v))\right)\\
&\hspace{-3mm}=&\hspace{-3mm}N\left(NT(u)T(v)+T(u)NT(v)-N(T(u)T(v))\right)\\
&\hspace{-3mm}\oeq{eq2.2}&\hspace{-3mm}NT(u)\cdot NT(v).
\end{eqnarray*}

(3) Note that $TS=NT$ and the previous statement proves that $NT\in\K(A_V)$, thus $TS$ belongs to $\K(A_V)$ as well. 
We would like to show that $T+TS\in \K(A_V)$. Clearly, 
$u\dm^{T+TS}v=u\dm^{T}v+u\dm^{TS}v=u\dm^{T}v+u\dm_{S}^{T}v$.  Thus
\begin{eqnarray*}
&&(T+TS)(u\dm^{T+TS}v)\\
&=&T(u\dm^{T}v)+T(u\dm_{S}^{T}v)+T(S(u\dm^{TS}v))+T(S(u\dm^{T}v))\\
&=&T(u\dm^{T}v)+T(S(u\dm^{TS}v))+T(S(u\dm^{T}v))+T\left(S(u)\dm^{T}v+u\dm^{T}S(v)-S(u\dm^{T}v)\right)\\
&=&T(u\dm^{T}v)+TS(u\dm^{TS}v)+T\left(S(u)\dm^{T}v+u\dm^{T}S(v)\right)\\
&\oeq{eq3.4}&T(u)T(v)+TS(u)TS(v)+TS(u)T(v)+T(u)TS(v)\\
&=&(T+TS)(u)(T+TS)(v).
\end{eqnarray*}
On the other hand, $(T+TS)(u\dm^{T+TS}v)\oeq{eq2.14}(T+TS)(\ell_{(T+TS)(u)}(v)+r_{(T+TS)(v)}(u))$. Hence,
$$(T+TS)(u)(T+TS)(v)=(T+TS)(\ell_{(T+TS)(u)}(v)+r_{(T+TS)(v)}(u))$$
which means that $T+TS$ is a Kupershmidt operator. Therefore, $T$ and $TS$ are compatible.

(4) We use the induction on $k$ to prove that $T_k\in \K(A_V)$. Clearly, $T_0=TS^0=T\in \K(A_V)$. Now assume that
$T_{k-1}\in \K(A_V)$. Then $T_k=TS^k=NTS^{k-1}=NT_{k-1}\in\K(A_V)$ by the second statement.

Since $T\in\K(A_V)$ and $TS=NT$, it follows that $TS^i=N^iT$ and
\begin{eqnarray*}
T(u\dm^{T}_{S^{i}}v)&=&T\left(S^{i}(u)\dm^{T}v+u\dm^{T}S^{i}(v)-S^{i}(u\dm^{T}v)\right)\\
&\oeq{eq3.4}&TS^{i}(u)T(v)+T(u)TS^{i}(v)-N^{i}T(u\dm^{T}v)\\
&=&N^{i}T(u)T(v)+T(u)N^{i}T(v)-N^{i}(T(u)T(v))\\
&=&T(u)\cdot_{N^{i}}T(v).
\end{eqnarray*}
By the first statement, $S\in\NN(V,\dm^T)$ and thus together (\ref{eq3.5}) with $\Delta=0$ in the proof of the first statement implies that $S(u\dm_{S}^{T}v)=S(u\dm^{TS}v)=S(u)\dm^{T}S(v).$ Using induction on $k\in\N^+$, we see that
\begin{equation}
\label{eq3.7}
S^k(u\dm_{S^k}^{T}v)=S^k(u)\dm^{T}S^k(v).
\end{equation}
In fact, note that $TS=NT\in\K(A_V)$ and we have seen that the case of $k=1$ holds. Thus by induction hypothesis, it follows that $S^k(u\dm_{S^k}^{T}v)=S\left(S^{k-1}(u\dm_{S^{k-1}}^{TS}v)\right)=
S\left(S^{k-1}(u)\dm^{TS}S^{k-1}(v)\right)=S\left(S^{k-1}(u)\dm^{T}_SS^{k-1}(v)\right)=S^k(u)\dm^{T}S^k(v)$.
Thus, (\ref{eq3.7}) holds. 

More generally, for all $i\in\N$, it follows that
\begin{equation}
\label{eq3.8}
S^k(u\dm_{S^{k+i}}^{T}v)=S^k(u\dm_{S^{k}}^{TS^i}v)\oeq{eq3.7}S^k(u)\dm^{TS^i}S^k(v)=S^k(u)\dm^{T}_{S^i}S^k(v).
\end{equation}
Hence,
$T_k(u\dm_{S^{k+i}}^{T}v) = T(S^k(u\dm_{S^{k+i}}^{T}v))  \oeq{eq3.8} T(S^{k}(u)\dm^{T}_{S^{i}}S^{k}(v))=TS^{k}(u)\cdot_{N^i}TS^{k}(v)=T_{k}(u)\cdot_{N^i}T_{k}(v).$

(5) We may suppose that $j=k+i$ for some $i\in\N$. We only need to show that $T_k+T_j=T_k+T_{k+i}\in\K(A_V)$. Note that $u\dm^{T_k+T_{k+i}}v=u\dm^{T_k}v+u\dm^{T_{k+i}}v=u\dm^{T_k}v+u\dm^{T_k}_{S^i}v$. Then
\begin{eqnarray*}
&&(T_k+T_{k+i})(u\dm^{T_k+T_{k+i}}v)\\
&=&(T_k+T_{k+i})(u\dm^{T_k}v+u\dm^{T_k}_{S^i}v)\\
&=&T_k(u\dm^{T_k}v)+T_{k}(u\dm^{T_k}_{S^i}v)+T_{k+i}(u\dm^{T_k}v)+T_{k+i}(u\dm^{T_k}_{S^i}v)\\
&\oeq{eq2.15}&T_k(u)T_k(v)+T_{k}\left(S^{i}(u)\dm^{T_k}v+u\dm^{T_k}S^{i}(v)-S^{i}(u\dm^{T_k}v)\right)+T_{k+i}(u\dm^{T_k}v)+T_{k+i}(u\dm^{T_{k+i}}v)\\
&=&T_k(u)T_k(v)+T_{k}\left(S^{i}(u)\dm^{T_k}v+u\dm^{T_k}S^{i}(v)-S^{i}(u\dm^{T_k}v)\right)+T_{k}(S^i(u\dm^{T_k}v)+T_{k+i}(u\dm^{T_{k+i}}v)\\
&\oeq{eq3.4}&T_k(u)T_k(v)+T_{k}(S^{i}(u))T_{k}(v)+T_{k}(u)T_{k}(S^{i}(v))+T_{k+i}(u)T_{k+i}(v)\\
&=&T_k(u)T_k(v)+T_{k+i}(u)T_{k}(v)+T_{k}(u)T_{k+i}(v)+T_{k+i}(u)T_{k+i}(v)\\
&=&(T_k+T_{k+i})(u)(T_k+T_{k+i})(v).
\end{eqnarray*}
On the other hand, it follows from (\ref{eq2.14}) that 
$(T_k+T_{k+i})(u)(T_k+T_{k+i})(v)=(T_k+T_{k+i})(u\dm^{T_k+T_{k+i}}v)=(T_k+T_{k+i})(\ell_{(T_k+T_{k+i})(u)}(v)+r_{(T_k+T_{k+i})(v)}(u)).$ Hence, $T_k+T_{k+i}\in\K(A_V)$.
\end{proof}

\section{More Examples and Computations}\label{sec4}
\setcounter{equation}{0}
\renewcommand{\theequation}
{4.\arabic{equation}}
\setcounter{theorem}{0}
\renewcommand{\thetheorem}
{4.\arabic{theorem}}

\noindent This section is devoted to providing more examples and detailed calculations in low dimensions.

\subsection{Low-dimensional pre-Malcev algebras}

In Example \ref{exam2.2}, we have seen a family $A(a,b)$ of two-dimensional nontrivial pre-Malcev algebras over $\C$.
A natural and interesting question might ask how to classify all two-dimensional nontrivial complex pre-Malcev algebras.
Compared with a classification of two-dimensional complex pre-Lie algebras (see for example, \cite[Theorem 3.1]{ZB12}), our computational experiences lead us to make the following conjecture. 

\begin{conj}\label{conj4.1}
{\rm
Let $A$ denote a nontrivial two-dimensional complex pre-Malcev algebra spanned by $\{e_1,e_2\}$. Then the multiplication of $A$ must be given by the following form:
$$e_1e_1=ae_1+c_1e_2, e_1e_2=c_2e_1-ae_2, e_2e_1=-be_1+c_3e_2,e_2e_2=c_4e_1+be_2$$
where $a,b,c_1,\dots,c_4\in\C$ satisfying 
\begin{eqnarray*}
 c_2^2 + 2bc_2-c_3c_4 &=&0\\
c_3^2+2ac_3- c_1c_2&=&0\\
c_1c_4  - c_2c_3 - 2ac_2 - 2bc_3 - 4ab&=&0.
\end{eqnarray*}
}\end{conj}

Looking back to \cite{Bai09} where all three-dimensional pre-Lie algebras over $\C$ were classified, it is reasonable to believe that classifying three-dimensional nontrivial complex pre-Malcev algebras would be more complicated. Here we present several families of three-dimensional nontrivial complex pre-Malcev algebras that might be helpful to the classification problem of low-dimensional nontrivial pre-Malcev algebras. 

In the following three examples, we suppose that $A$ is a three-dimensional complex nonassociative algebra spanned by $\{e_1,e_2,e_3\}$. These examples can be verified directly using the definitions of pre-Malcev algebras and pre-Lie algebras.

\begin{exam} \label{exam4.2}
{\rm
The following nonzero multiplication relations
$$e_1e_1=e_2+e_3, e_3e_3=ae_2$$
where $a\in\C^\times$, define a nontrivial pre-Malcev algebra structure on $A$.
\hbo}\end{exam}

\begin{exam} \label{exam4.3}
{\rm
The following nonzero multiplication relations
$$e_1e_1=ae_1, e_1e_2=-ae_2, e_2e_2=be_1, e_3e_3=ce_3$$
where $a,b,c\in\C^\times$, also define a nontrivial pre-Malcev algebra structure on $A$.
\hbo}\end{exam}

\begin{exam} \label{exam4.4}
{\rm
Suppose that all nonzero multiplications in $A$ are given by
$$e_1e_1=ae_1+be_3, e_1e_2=-ae_2+ce_3, e_2e_2=de_1+\frac{bd}{a}e_3$$
for some $a,b,c,d\in\C^\times$. Then $A$ is a nontrivial pre-Malcev algebra.
\hbo}\end{exam}

A four-dimensional example of nontrivial pre-Malcev algebras has been studied in \cite{RZ24}. 
Here we construct a family of four-dimensional pre-Malcev algebras that contains this example. 
A generalization of this example in theory of Hom-pre-Malcev algebras can be found in 
\cite[Example 3.1]{HMNS21}.

\begin{exam} \label{exam4.5}
{\rm
Let $A$ be a $4$-dimensional complex nonassociative  algebra spanned by $\{e_1,e_2,e_3,e_4\}$ with the following 
nonzero multiplication relations:
$$e_1e_2=e_2, e_1e_3=e_3,e_1e_4=-e_4, e_2e_2=ae_4, e_2e_3=2e_4,e_3e_3=be_2$$
where $a,b\in\C$. Then $A$ is a nontrivial pre-Malcev algebra over $\C$. In particular, when $(a,b)=(0,0)$, this example has been studied in \cite[Example 4.1]{RZ24}.
\hbo}\end{exam}

\subsection{Kupershmidt operators  and Nijenhuis operators}
In this subsection, we take $A$ to be the three-dimensional pre-Malcev algebra in Example \ref{exam4.3} with
$a=b=c=1$. In other words, the nonzero multiplications in $A$ are given by
$$e_1e_1=e_1, e_1e_2=-e_2, e_2e_2=e_1, e_3e_3=e_3.$$
This subsection aims to compute the variety $\K(A)$ of Rota-Baxter operators (i.e., Kupershmidt operators associated to the adjoint bimodule) and the variety $\NN(A)$ of Nijenhuis operators on $A$, providing a geometric understanding of these two varieties. 

We will give all details in describing the geometric structure of $\K(A)$ in Corollary \ref{coro4.11} and apply a similar method to characterize the geometry of $\NN(A)$. Thus we only  state the main result about the geometric structure of $\NN(A)$ in Theorem \ref{thm2} below, without going into detailed proofs. 

Let us focus on $\K(A)$. Suppose that 
$$T=\begin{pmatrix}
    a_{11}  & a_{12} & a_{13}  \\
     a_{21}  & a_{22} & a_{23}  \\
     a_{31}  & a_{32} & a_{33}  \\  
\end{pmatrix}\in\K(A)$$
denotes a generic Kupershmidt operator on $A$ and the action $T$ on $A$ is given by
\begin{eqnarray*}
T(e_1) & = & a_{11}e_1+  a_{12}e_2 + a_{13}e_3 \\
T(e_2) & = & a_{21}e_1+  a_{22}e_2 + a_{23}e_3 \\
T(e_3) & = & a_{31}e_1+  a_{32}e_2 + a_{33}e_3.
\end{eqnarray*}
Since $T$ satisfies (\ref{eq2.18}) with $(x,y)$ running over $\{(e_i,e_j)\mid 1\leqslant i,j\leqslant 3\}$, this means that
$T$ can be understood as a point in the affine space $\C^9$. Thus, $\K(A)$ can be viewed as an affine subvariety 
of $\C^9$ and the vanishing ideal $J$ of $\K(A)$ will be a radical ideal of 
$\C[x_{ij}\mid 1\leqslant i,j\leqslant 3]$.

To determine the ideal $J$, we need to compute (\ref{eq2.18}) when $(x,y)$ runs over $\{(e_i,e_j)\mid 1\leqslant i,j\leqslant 3\}$ and we will obtain 25 polynomial equations in $a_{ij}$ after deleting redundant ones. We write 
$I$ for the ideal generated by the corresponding 25 polynomials in $\C[x_{ij}\mid 1\leqslant i,j\leqslant 3]$. 
Then $J=\sqrt{I}$ and $V(J)=V(I)$. Define
\begin{eqnarray*}
 f_1:=    x_{11} + x_{22} & f_2:= x_{12} - 2x_{21}&
  f_3:=  x_{13}x_{22} + 2x_{23}^2\\
  f_4:= x_{13}x_{23} - x_{22}x_{23} &
    f_5:=   x_{13}^2 + 2x_{23}^2   & f_6:= x_{13}x_{21} + x_{22}x_{23}\\
    f_7:=  2x_{21}^2 + x_{22}^2 & f_8:=x_{21}x_{23} + x_{23}^2
    &f_9:=x_{22}^2x_{23} + 2x_{23}^3.
\end{eqnarray*}
Then $J=\ideal{x_{3i},f_j\mid 1\leqslant i\leqslant 3,1\leqslant j\leqslant 9}$ and all evaluations of all generators of $J$ on
$\K(A)$ are zero because $J$ is the radical of $I$ and the generators of $I$ come from the 25 polynomial equations. Thus $\K(A)\subseteq V(J).$ We would like to show that
\begin{equation}
\label{ }
\K(A)=V(J).
\end{equation}

Define $\p_1:=\ideal{x_{31},x_{32},x_{33}, x_{13},x_{23}, f_1, f_2, f_7}$
and $$\p_2:=\ideal{x_{31},x_{32},x_{33}, f_1, g:=x_{12} + 2x_{23},
f_4/x_{23}, f_8/x_{23}, f_9/x_{23}}.$$

\begin{lem}\label{lem4.6}
The ideals $\p_1$ and $\p_2$ are prime.
\end{lem}

\begin{proof}
To show that an ideal is prime, it suffices to show that the corresponding quotient ring is 
an integral domain. Note that 
$$\C[x_{ij}\mid 1\leqslant i,j\leqslant 3]/\p_1\cong \C[x_{11},x_{12},x_{21},x_{22}]/\ideal{f_1, f_2, f_7}\cong\C[x_{21},x_{22}]/\ideal{f_7}$$
which is an integral domain because $f_7$ is irreducible in $\C[x_{21},x_{22}]$. This implies that $\p_1$ is a prime ideal. 
Moreover, 
$$\C[x_{ij}\mid 1\leqslant i,j\leqslant 3]/\p_2\cong \C[x_{11},x_{12},x_{13},x_{21},x_{22},x_{23}]/\ideal{f_1, g,
f_4/x_{23}, f_8/x_{23}, f_9/x_{23}}$$
which is isomorphic to $\C[x_{22},x_{23}]/\ideal{x_{22}^2+ 2x_{23}^2}$ that is an integral domain as well. Hence,
$\p_2$ is a prime ideal.
\end{proof}

\begin{lem}\label{lem4.7}
$V(\p_1)\cup V(\p_2)\subseteq V(J)$.
\end{lem}

\begin{proof}
It suffices to show that $J\subseteq \p_1$ and $J\subseteq \p_2$. To prove that the first containment, we need to verify that
each element in $\{f_i\mid 1\leqslant i\leqslant 9\}\setminus\{f_1,f_2,f_7\}$ is equal to zero modulo $\p_1$. This is immediate because $f_4,f_8,f_9$ are divisible by $x_{23}$ and each term of anyone of $\{f_3,f_5,f_6\}$ is divisible by either 
$x_{23}$ or $x_{13}$, while $x_{13}\equiv 0\equiv x_{23}$ modulo $\p_1$. Moreover, 
to show that $J\subseteq \p_2$, it suffices to show that $f_2,f_3,f_5,f_6,f_7$ are equal to zero modulo $\p_2$. We may be working over modulo $\p_2$. Note that $g\in\p_2$, thus $x_12\equiv -2x_{23}$ and $f_2=-2(x_{21}+x_{23})=-2(f_8/x_{23})\equiv 0$ because $f_8/x_{23}$ lies in $\p_2$. The fact that $f_4/x_{23}\in \p_2$ implies that
$x_{13}\equiv x_{22}$, thus $f_3\equiv f_5\equiv x_{22}^2+2x_{23}^2=f_9/x_{23}\equiv 0$.  The fact that $f_8/x_{23}\in \p_2$ implies that $x_{21}\equiv -x_{23}$, thus $f_6\equiv -x_{22}x_{23}+x_{22}x_{23}=0.$ For $f_7$, we see that
it is equivalent to $x_{22}^2+2x_{23}^2=f_9/x_{23}\equiv 0$. This proves that $J$ is also contained in $\p_2$.
\end{proof}

\begin{lem}\label{lem4.8}
$V(J)\subseteq V(\p_1)\cup V(\p_2)$.
\end{lem}

\begin{proof}
Note that by \cite[Theorem 7, page 192]{CLO15}, we see that $V(\p_1)\cup V(\p_2)=V(\p_1\cdot\p_2)$. Thus it suffices to show that $\p_1\cdot\p_2\subseteq J$. It follows from \cite[Proposition 6, page 191]{CLO15} that $\p_1\cdot\p_2$ can be generated by
the products of the generators for $\p_1$ and $\p_2$. Thus we only need to show that
$$\{ab\mid a\in\{x_{13},x_{23}\},b\in \{g,
f_4/x_{23}, f_8/x_{23}, f_9/x_{23}\}\}\subseteq J.$$
Let us work over modulo $J$. Note that $f_2,f_4,f_6\in J$, thus
$x_{12} \equiv 2x_{21}, x_{13}x_{23}\equiv x_{22}x_{23}$, and $x_{13}x_{21}\equiv -x_{22}x_{23}$.
Hence, $x_{13}g=x_{13}(x_{12} + 2x_{23})=x_{12} x_{13}+ 2x_{13}x_{23}\equiv 2x_{21} x_{13}+ 2x_{22}x_{23}=-2x_{22}x_{23}+2x_{22}x_{23}=0.$ Namely, $x_{13}g\in J$. Similar arguments show that other products belong to $J$ as well.
\end{proof}

Together Lemmas \ref{lem4.6}, \ref{lem4.7}, and \ref{lem4.8} obtain the following results.

\begin{coro}\label{coro4.9}
$V(J)=V(\p_1)\cup V(\p_2)$ is an irreducible decomposition of affine varieties. 
\end{coro}

In particular, the following corollary shows the shapes of elements in $V(\p_1)$ and $V(\p_2)$.

\begin{coro}\label{coro4.10}
An element $T_1\in V(\p_1)$ is of the form
$$T_1=\begin{pmatrix}
    a & 2b & 0  \\
     b  & -a & 0  \\
     0  & 0 & 0  \\  
\end{pmatrix}$$
with $a^2+2b^2=0$ and an element $T_2\in V(\p_2)$ has the form
$$T_2=\begin{pmatrix}
    a & -2b & a  \\
    -b   & -a & b  \\
     0  & 0 & 0  \\  
\end{pmatrix}$$
where $a^2+2b^2=0$.
\end{coro}

\begin{coro}\label{coro4.11}
$\K(A)=V(\p_1)\cup V(\p_2)$. In other words, $\K(A)$ is a one-dimensional affine variety consisting of two irreducible components.
\end{coro}

\begin{proof}
Since $\K(A)\subseteq V(J)$, it follows from Corollary \ref{coro4.9} that $\K(A)\subseteq V(\p_1)\cup V(\p_2)$. 
A direct verification shows that $T_1$ and $T_2$ in Corollary \ref{coro4.10} both are Kupershmidt operators on $A$. Thus
$V(\p_1)\cup V(\p_2)\subseteq \K(A)$.
\end{proof}

Applying the same method, we derive the following result that characterizes 
 the geometric structure of $\NN(A)$ without going into detailed proofs.

\begin{thm}\label{thm2}
The variety $\NN(A)$ is a $3$-dimensional affine variety consisting of $7$ irreducible components $\{V(\p_i)\mid 1\leqslant i\leqslant 7\}$ for which $\dim(V(\p_1))=\dim(V(\p_2))=3$ and $\dim(V(\p_i))=2$ for $3\leqslant i\leqslant 5$, where 
generic elements $N_i\in V(\p_i)$ have the following forms:
\begin{eqnarray*}
N_1&=&\begin{pmatrix}
    a & 0 & 0 \\
     b & c & 0  \\
     0  & 0 & d \\  
\end{pmatrix}\textrm{ with }a^2 - 2ac + 2b^2 + c^2\\
N_2&=&\begin{pmatrix}
    a & -2b & 0 \\
     b  & a & 0  \\
    0 & 0& c \\  
\end{pmatrix}\\
N_3&=&\begin{pmatrix}
    a+c  & 0& c  \\
     b & a & b  \\
     0 & 0 & a  \\  
\end{pmatrix} \textrm{ with }c^2+2b^2=0\\
N_4&=&\begin{pmatrix}
    a & 2b & c-a \\
     -b  & a& b  \\
     0 & 0 & c  \\  
\end{pmatrix} \textrm{ with }a^2 - 2ac + 2b^2 + c^2=0\\
N_5&=& \begin{pmatrix}
    a-b & 0& 0 \\
     0  & a-b& 0  \\
     b & 0 & a  \\  
\end{pmatrix}\\
N_6&=& \begin{pmatrix}
    a+c  & 0 & 0  \\
     0  & a+c& 0  \\
     a  & b& c \\  
\end{pmatrix} \textrm{ with }2a^2+ b^2=0\\
N_7&=& \begin{pmatrix}
    a  & 0 & 0  \\
     -b & c & b  \\
     0  & 0 & a \\  
\end{pmatrix} \textrm{ with }c^2 - 2ac + 2b^2 + a^2=0
\end{eqnarray*}
for $a,b,c,d\in \C$.
\end{thm}

We close this article with the following remark, which was suggested by the referee’s comments.

\begin{rem}{\rm
Note that all pre-Malcev algebras are  Malcev-admissible algebras. Thus it is natural to ask whether a Nijenhuis 
(or Kupershmidt) operator on a pre-Malcev algebra $A$ is still a Nijenhuis 
(or Kupershmidt) operator on the corresponding Malcev algebra $[A]$. In fact, if $A$ is a pre-Lie algebra, then $[A]$ is a Lie algebra. By \cite[Proposition 4.7]{WSBL19}, we have seen that a Nijenhuis operator on $A$ is also a Nijenhuis operator on 
$[A]$. This result can be generalized to the case of pre-Malcev algebras via a direct verification, i.e.,  Nijenhuis operators on 
a pre-Malcev algebra $A$ are Nijenhuis operators on the corresponding Malcev algebra $[A]$. However, the converse 
statement might not be true in general. For example, let us consider Example \ref{exam2.4}. The corresponding Malcev algebra
$[A]$ is defined by the nontrivial relations $[e_1,e_2]=-e_2$ and $[e_2,e_1]=e_2$. Direct computations show that
any $2$ by $2$ matrix is a Nijenhuis operator on $[A(1,0)]$, while in Example \ref{exam2.4}, we see that not all matrices are  
Nijenhuis operators on $A(1,0)$. Similar conclusions can be obtained for Kupershmidt operators and Rota-Baxter operators;
see \cite[Corollary 2.12]{AB08} for more details.
\hbo}\end{rem}

\vspace{3mm}
\noindent \textbf{Acknowledgements}. 
The first-named author would like to thank Professors H. Eddy A. Campbell and David L. Wehlau for their helpful conversations and encouragement. This research was partially supported by NNSF of China (Grant No. 12561003). The authors would like to thank the two anonymous referees and the editor for their careful reading, constructive comments, and suggestions.

%%%%%%%%%%%%%%%%%%%%%%%%%%References%%%%%%%%%%%%%%%%%%%%%%%%

\raggedright
\end{document}